%% file: 0_main.tex
\documentclass[12pt]{article}
\usepackage{graphicx} % Required for inserting images

\title{\huge A Problem of Calculus of Variations and Game Theory}
\author{Christopher Boyer, Grace Luo, and Siddharth Penmetsa \footnote{\scriptsize{Mathematics Department, North Carolina School of Science and Mathematics (NCSSM), Durham, NC 27705, USA}}}
\date{}
\usepackage[letterpaper,top=2cm,bottom=2cm,left=3cm,right=3cm,marginparwidth=1.75cm]{geometry}
\usepackage[thinc]{esdiff}
\usepackage{amsmath}
\usepackage{amsthm}
\usepackage{float}
\usepackage{mathtools}
\usepackage{array}
\usepackage{blkarray}

\usepackage{biblatex}
\usepackage{hyperref}
\hypersetup{
    colorlinks=true,
    linkcolor=black,
    urlcolor=blue,
    citecolor=blue,
    }
\begin{document}

\maketitle

\begin{abstract}
In this paper, we study a theoretical math problem of game theory and calculus of variations in which we minimize a functional involving two players. A general relationship between the optimal strategies for both players is presented, followed by computer analysis as well as polynomial approximation. Nash equilibrium strategies are determined through algebraic manipulation and linear programming. Lastly, a variation of the game is also investigated.
\end{abstract}

\input{1_Introduction}
\input{2_CalculusOfVariations}
\input{3_ComputationalAnalysis}
\input{4_PolynomialApproximation}
\input{5_TwoLemmas}
\input{6_SmallA}
\input{7_LargeA}
\input{8_sin_3x}
\input{9_FutureWork}
\input{10_Acknowledgements}

\printbibliography

\end{document}

%% file: 1_Introduction.tex
\section{Introduction}
Our problem combines aspects of calculus of variations and game theory. The game operates as follows:
Define the functional $S$ such that $$S
\big(f(t), g(t)\big)=\int_{0}^{1}\bigg(f'(t)^2-g'(t)^2-a\sin\big(f(t)-g(t)\big)\bigg)\ dt,$$
where $a$ is a non-negative real number.\\

There are two players in the game, $f$ and $g$. Player $f$ will choose a differentiable function $f(t)$ and player $g$ will choose a differentiable function $g(t)$. There are a few conditions on $f(t)$ and $g(t)$:

$$f(0)=f'(1)=0$$
$$g(0)=g'(1)=0.$$

The game states that $S\big(f(t), g(t)\big)=S$ is the amount of money that $f$ will pay $g$. Therefore, $f$ wants to minimize $S$ and $g$ wants to maximize $S$. Since we want the game to be symmetric ($f$ playing against $g$ should yield the same exchange of money as if $g$ were playing against $f$), $S$ must be an odd function, such that $S\big(f(t), g(t)\big)=-S\big(g(t), f(t)\big)$. After all, if $f$ pays $g$ $k$ dollars, then $g$ is paying $-k$ dollars to $f$. Note that we could substitute the $\sin(t)$ function in the functional $S$ with any other odd function, such as $\tan^{-1}(t)$, and this condition would still be true. Such changes would result in variations of the game.\\

We seek to better understand the functional in this game, as well as optimal strategies for both players, and investigate if a Nash equilibrium exists. A Nash equilibrium is defined as a position of no regret. In other words, $f$ and $g$ would choose the same functions regardless of if they knew each other's moves or not \cite{nash}.

%% file: 2_CalculusOfVariations.tex
\section{Calculus of Variations}

\textbf{Problem}: Given a fixed $g(t)$,
find the function $f(t)$, for $0\le t\le1$, that satisfies the conditions $f(0)=f'(1)=0$ and minimizes $$S\bigl(f(t),\ g(t)\bigl)=S=\int_{0}^{1}\bigg(f'(t)^{2}-g'(t)^{2}-a\sin\big(f(t)-g(t)\big)\bigg)\ dt.$$

Let $f(t)=F(t)+s\eta(t)$ be a variation of $F(t)$, where $F(t)$ is the optimal solution and $\eta(t)\neq0$ \cite{calcofvariations}. We can then say that $f'(t)=F'(t)+s\eta'(t)$ and $\dot{f}(t)\coloneqq\frac{\partial}{\partial s}f(t)=\eta(t)$. Since $f(0)=f'(1)=0$, and the optimal function $F(t)$ also satisfies the conditions $F(0)=F'(1)=0$, the same must be true for $\eta(t)$: $\eta(0)=\eta'(1)=0$.\\

We can then substitute $f(t)=F(t)+s\eta(t)$ and $f'(t)=F'(t)+s\eta'(t)$ into $S$ to obtain $$S=\int_{0}^{1}\bigg(\big(F'(t)+s\eta'(t)\big)^2-g'(t)^2-a\sin\big(F(t)+s\eta(t)-g(t)\big)\bigg)\ dt=S(s).$$
Following integration, $S$ will be in terms of $t$ and $s$. After $t=0$ and $t=1$ are substituted in, $S$ will solely be in terms of $s$.\\

In order to find where $S$ obtains a minimum, we first must find where $S$ has a critical point. $S$ has a critical point when $\frac{\partial}{\partial s}S=0$, which occurs when $s=0$ since then $f(t)=F(t)+0\cdot\eta(t)=F(t)$. This means that $f(t)$ is equal to the optimal solution $F(t)$ when $s=0$, or when there is no variation.\\

Therefore, the following holds true:
$$0=\diff{}{s}S\bigg\rvert_{s=0}=\dot{S}\bigg\rvert_{s=0}=\int_{0}^{1}\frac{\partial}{\partial s}\biggl({f'(t)^2-g'(t)^2-a\sin\bigl(f(t)-g(t)\bigl)}\biggl)\bigg\rvert_{s=0}dt.$$
Notice that the derivative turned into a partial derivative once we moved it inside the integral because the expression inside is still in terms of $t$.\\
$$0=\dot{S}\bigg\rvert_{s=0}=\int_{0}^{1}\biggl(2f'(t)\dot{f}'(t)-a\cos\bigl(f(t)-g(t)\bigl)\dot{f}(t)\biggl)\bigg\rvert_{s=0}dt$$\\

Using integration by parts, we find that $$\int_{0}^{1}f'(t)\dot{f}'(t)\ dt = f'(t)\dot{f}(t)\bigg\rvert_{0}^{1}-\int_{0}^{1}f''(t)\dot{f}(t)\ dt.$$
But, recall that $f'(1)=0$ is an initial condition, and since $\eta(0)=0$, $\dot{f}(0)=\eta(0)=0$. Therefore, $f'(t)\dot{f}(t)\bigg\rvert_{0}^{1}=0$, and $$\int_{0}^{1}f'(t)\dot{f}'(t)\ dt = \int_{0}^{1}-f''(t)\dot{f}(t)\ dt.$$

Substituting this in, we determine that $$0=\int_{0}^{1}\bigg(-2f''(t)\dot{f}(t)-a\cos(f(t)-g(t))\dot{f}(t)\bigg)\bigg\rvert_{s=0}dt \qquad \qquad \forall \dot{f}(t)$$

$$0=\int_{0}^{1}\bigg(-2f''(t)-a\cos(f(t)-g(t))\bigg)\dot{f}(t)\bigg\rvert_{s=0}dt \qquad \qquad \forall \dot{f}(t).$$

Recall that $f(t)=F(t)+s\eta(t)$; therefore, evaluating the inside of the integral at $s=0$, we obtain

$$0=\int_{0}^{1}\bigg(-2F''(t)-a\cos(F(t)-g(t))\bigg)\eta(t)\ dt \qquad \qquad \forall \eta(t).$$

The equation above must be true for any $\eta(t)$, meaning $$-2F''(t)-a\cos(F(t)-g(t))=0.$$

We can rewrite the equation above to determine the following second-order differential equation regarding $F(t)$, which is the optimal solution:
$$F''(t)=-\frac{a}{2}\cos\bigl(F(t)-g(t)\bigl).$$

However, since $f(t)=F(t)$ when $s=0$, we will refer to the optimal solution as $f(t)$ in the rest of the paper.

%% file: 3_ComputationalAnalysis.tex
\section{Computational Analysis}
The following graphs were created using Python and Desmos to explore how the second order differential equation behaves, which will allow us to predict solutions to the problem. Additionally, the graphs demonstrate the shooting method of finding solutions. We started with a list of $f(t)$ values satisfying $t=1$ and $f'(t)=0$ \cite{shootingmethod}. We then applied Euler's method with an extremely small step-size of $dt = 0.000001$. As these points move back in time (i.e., as $t$ approaches 0), we collected data on their positions (their $f(t)$ and $f'(t)$ values). \cite{bvp}\\

For simplicity, $g(t)$ was set to 0. We graphed the second order differential equation $f''(t)=-\frac{a}{2}\cos(f(t))$ in the phase plane, which made the solutions easier to visualize \cite{phaseplane}. In the graphs below, there are three axes. The red axis is time (or $t$), the green axis is $f(t)$, and the black axis is $f'(t)$.

\subsection{When $a=1$}
When $a=1$, the second order differential equation is $f''(t)=-\frac{a}{2}\cos(f(t))=-\frac{1}{2}\cos(f(t)).$

\begin{figure}[H]
    \centering
    \includegraphics[width=0.9\linewidth]{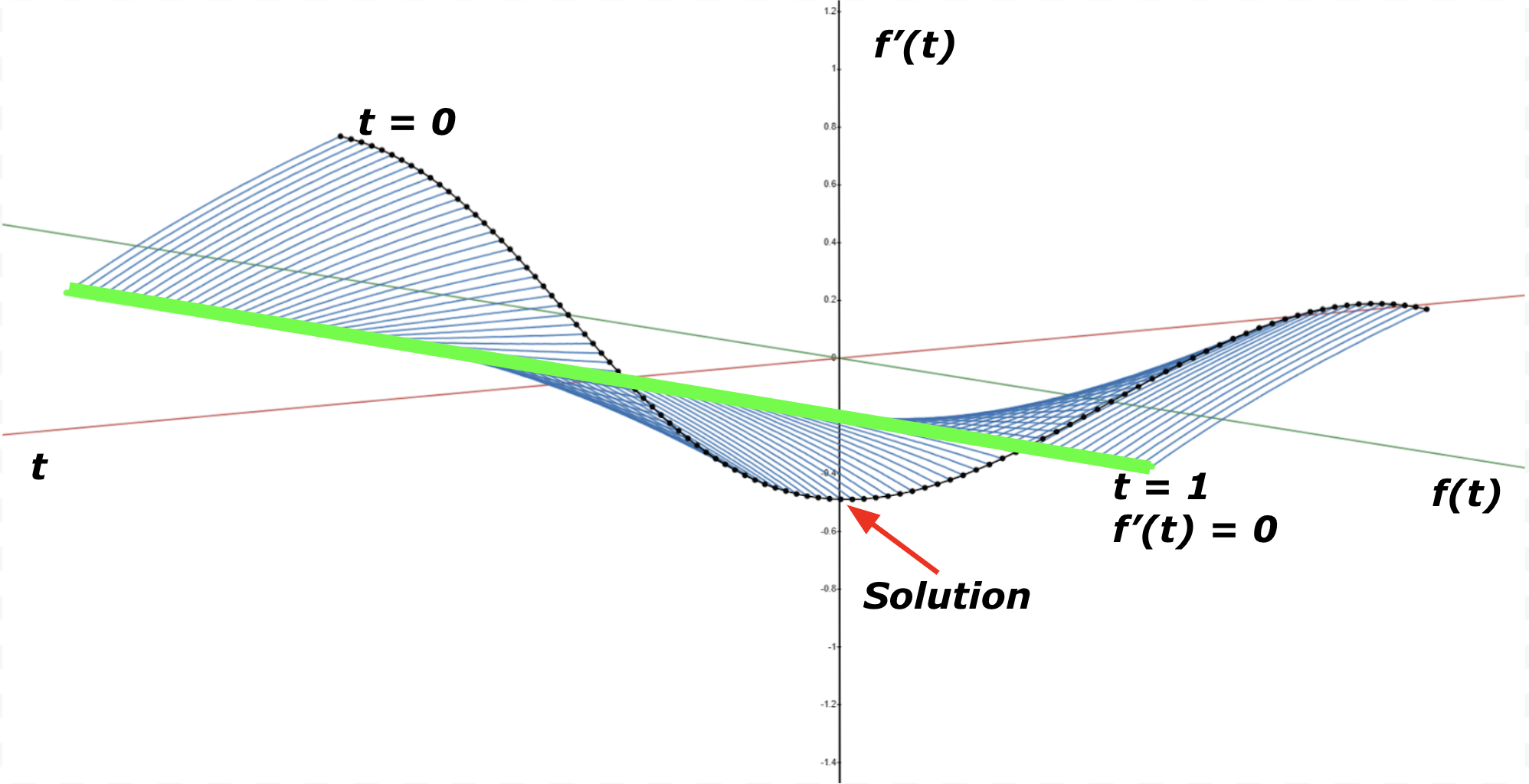}
    \caption{Labeled stream plot for $a=1$, where the neon line represents when $t=1$ and $f'(t)=0$, and the black dots represent $f(t)$ when $t=0$.}
    \label{fig:Figure 2}
\end{figure}

Every blue stream line in Figure 1 starts from the neon-green line where $t=1$ and $f'(t)=0$, which is one of the initial conditions, and stops when $t=0$, which is shown by the collection of black dots. The points on the neon-green line have $f(t)$ values that range from $-\pi$ to $\pi$ to show one complete cycle of radians. This trivial example illustrates how the shooting method finds solutions to the differential equation.
\begin{figure}[H]
    \centering
    \includegraphics[width=0.9\linewidth]{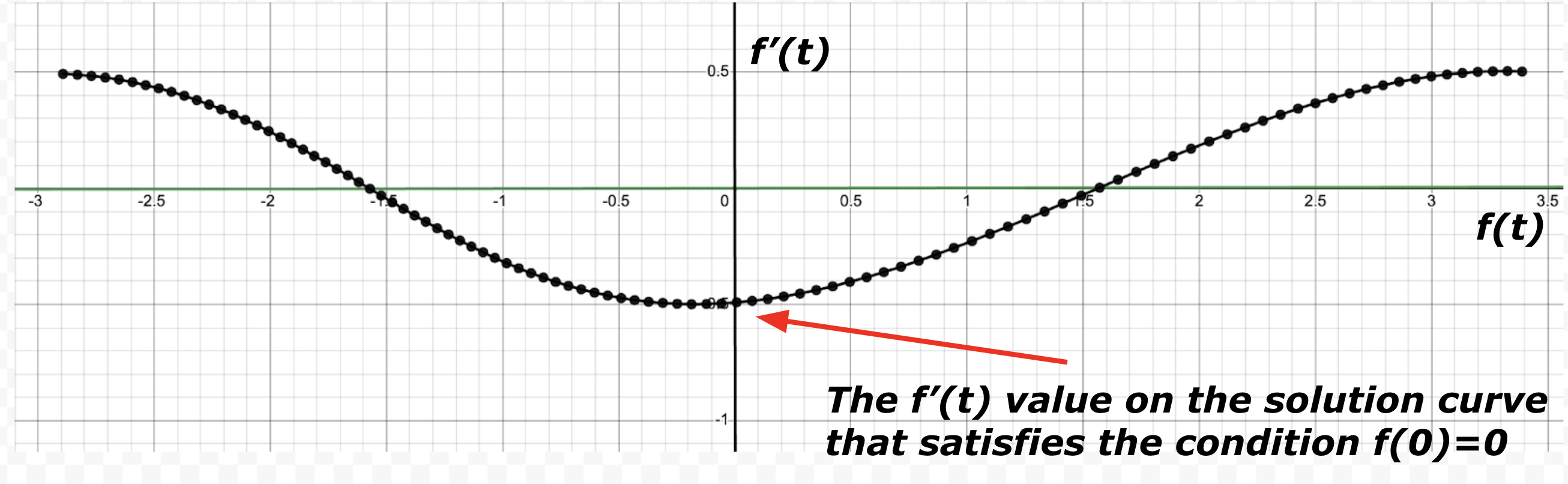}
    \caption{Illustrates all $f(t)$ and $f'(t)$ values for $a=1$ and $t=0$.}
    \label{fig:Figure 3}
\end{figure}

As shown in Figure 2, in the trivial case of $a=1$, there is only one solution that satisfies the initial condition $f(0)=0$, and this is when the solution curve intersects the $f'(t)$ axis (the vertical black axis).

\subsection{When $a=102$}
When $a=102$, the differential equation becomes $f''(t)=-\frac{a}{2}\cos(f(t))=-\frac{102}{2}\cos(f(t))=-51\cos(f(t))$.\\

Note: In the next few graphs, $f'(t)$ is scaled down by a factor of 10 for easier readability.

\begin{figure}[h]
    \centering
    \includegraphics[width=0.9\linewidth]{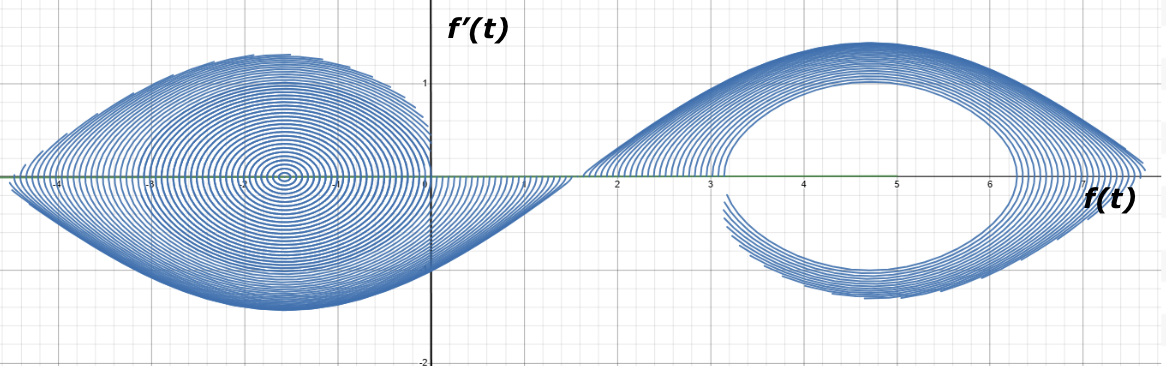}
    \caption{Stream plot of $f(t)$ from $t=0$ to $t=1$ for $a=102$.}
    \label{fig:Figure 4}
\end{figure}

Like before, the blue stream lines are potential paths of $f(t)$ from $t=0$ to $t=1$. The stream plot shows the washing machine effect from the differential equation \cite{phasegraph}. This demonstrates that multiple solutions exist, since multiple stream lines intersect the $f'(t)$ axis (which is when $f(0)=0$), and it is necessary to consider them all in game-play.

\begin{figure}[H]
    \centering
    \includegraphics[width=0.9\linewidth]{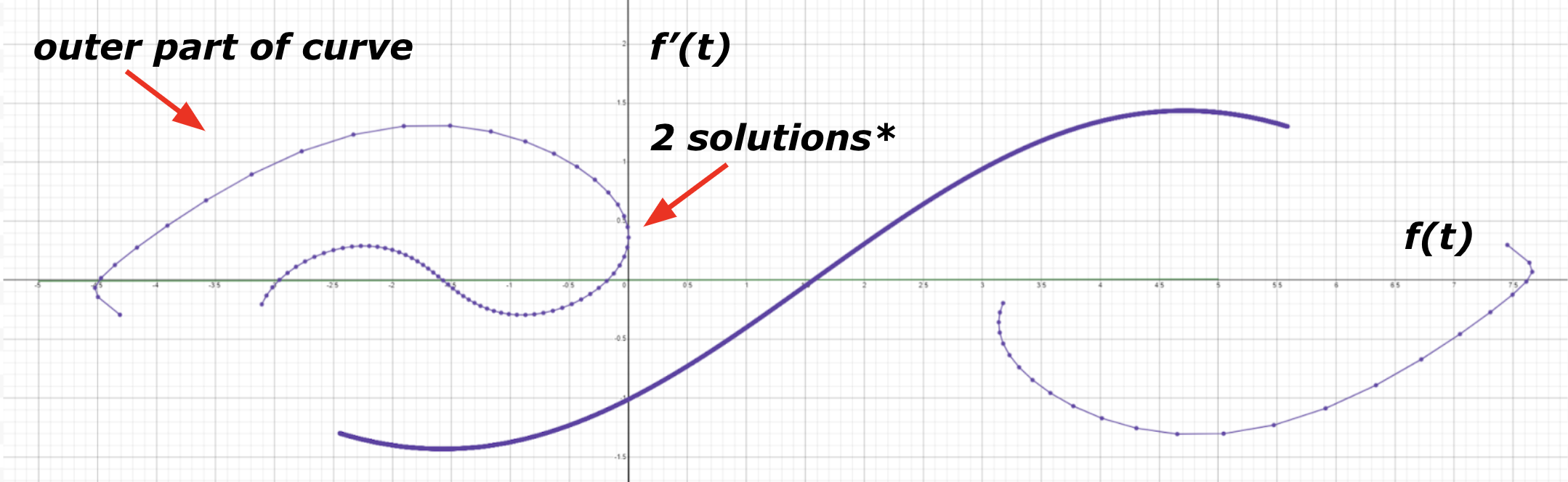}
    \caption{$f(t)$ when $t=0$ for the stream plot where $a=102$, with three solutions that satisfy the initial conditions. (*It may look like the purple curve is tangent to the $f'(t)$ axis, but if you zoom in, there are actually two solutions there.)}
    \label{fig:Figure 5}
\end{figure}

In Figure 4, we simplified the stream plot so that only the ``startpoints,” or the points when $t=0$, are shown by the broken up purple line. The purple line shows all possible starting points, but the critical points, or the solutions to the differential equation, occur only when the function intersects the vertical axis, which is the $f'(t)$ axis. (Note: We graphed less data points on the outer two curves due to Desmos’s list size limit.)\\

The broken up purple lines show that there are 3 solutions for $a=102$ when $f(t)=0$ at $t=0$. This strengthens our understanding of the differential equation and allows us to anticipate these solutions when solving the game theory problem.

\subsection{When $a=1000$}
When $a=1000$, the differential equation becomes $f''(t)=-\frac{a}{2}\cos(f(t))=-\frac{1000}{2}\cos(f(t))=-500\cos(f(t))$.\\

\begin{figure}[h]
    \centering
    \includegraphics[width=0.9\linewidth]{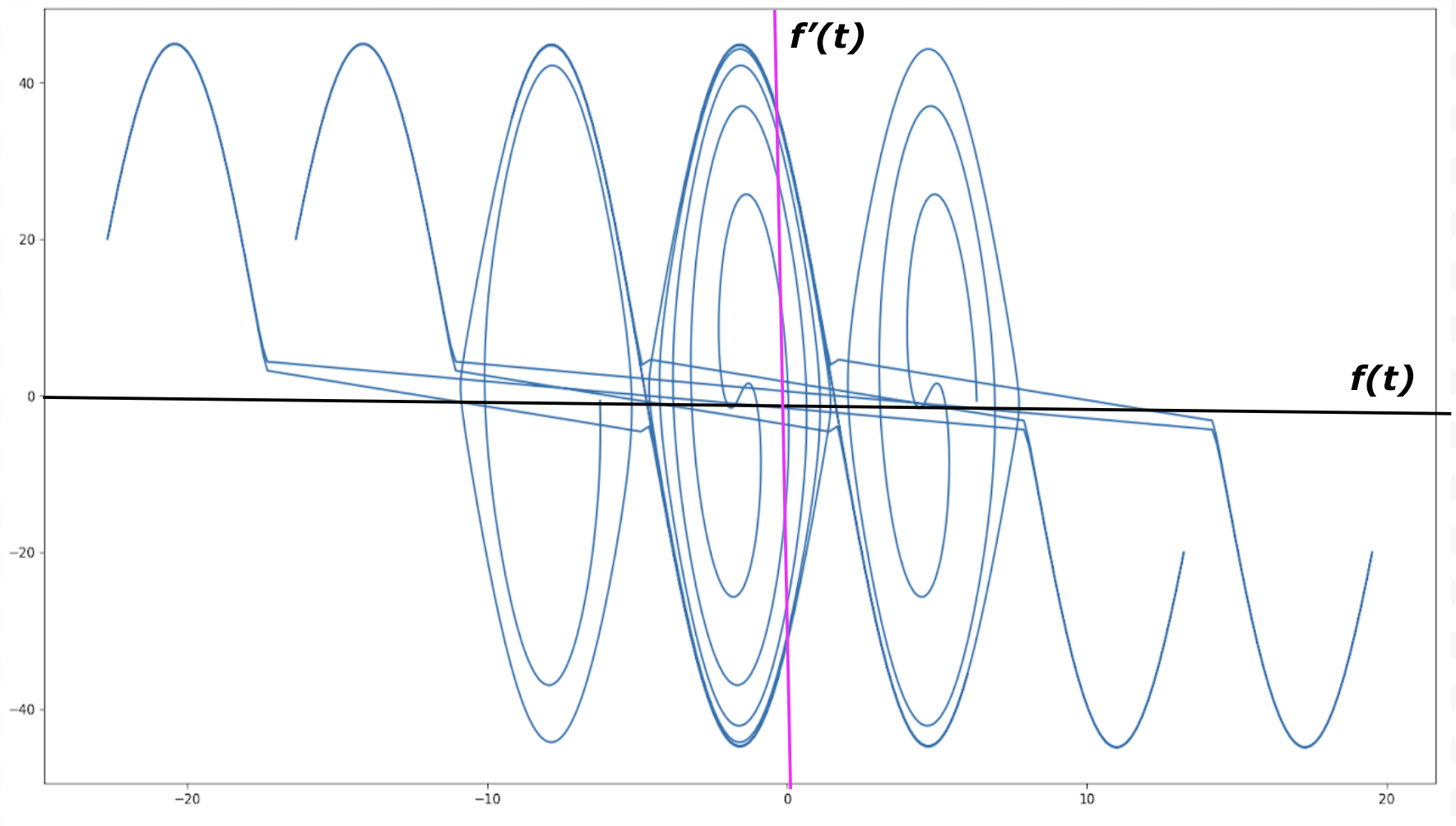}
    \caption{All possible ``startpoints" (when $t=0$) for the case where $a=1000$ (the straight lines in the middle of the graph are not supposed to be there, this was just how the datapoints were generated, and they show the jump in the $(f(t),f'(t))$ values at $t=0$).}
    \label{fig:Figure 6}
\end{figure}

There are around 12-14 solutions for $a=1000$ and potentially more, with endpoints (when $t=1$) ranging from $f(t)=-2\pi$ to $2\pi$.\\

These plots illustrate how increasing the $a$ value of the differential equation increases the number of solutions. We can make the conjecture that as $a\rightarrow\infty$, so does the number of solutions to the differential equation. With high values of $a,$ the fact that there are multiple solutions suggests that there are going to be multiple complex and mixed strategies for playing the game, and that the game could potentially turn into something that mimics a weighed rock-paper-scissors game, where the Nash equilibrium involves three functions, instead of just one. \cite{glico}

%% file: 4_PolynomialApproximation.tex
\section{Polynomial Approximation}
From calculus of variations, we found that the solution to the functional $S$ for any $a$ value has the second-order differential equation $$f''(t)=-\frac{a}{2}\cos\bigl(f(t)-g(t)\bigl).$$

For simplicity, let us assume that $g(t)=0$. Then, we can approximate $f''(t)$ by using 
the Taylor Series

$$f''(t)=-\frac{a}{2}\cos\bigl(f(t)\bigl)\ \approx -\frac{a}{2}\bigl(1-\frac{f(t)^2}{2}\bigl),$$

assuming that the rest of the higher order terms will be negligible since the value of $a$ is correlated with the value of $f(t)$, and is taken to an increasingly high power.\\

Now, we will attempt to model $f(t)$ using a polynomial approximation to better understand the behavior of the function. It is reasonable to claim that there is a unique polynomial, with up to an infinite number of terms, and thus an infinite number of degrees of freedom, that can accurately model $f(t)$.\\

Let $f(t)=b_0+b_1t+b_2t^2+b_3t^3+b_4t^4+b_5t^5+...$\\

Then, we can write the following:

$$f'(t)=b_1+2b_2t+3b_3t^2+4b_4t^3+5b_5t^4+...$$

$$f''(t)=2b_2+6b_3t+12b_4t^2+20b_5t^3+...$$

$$f(t)^2=b_0^2+(2b_0b_1)t+(2b_0b_2+b_1^2)t^2+(2b_0b_3+2b_1b_2)t^3+...$$\\

We can also use the initial condition $f(0)=0$ to give us more information regarding the coefficients in the polynomial. We know that $f(0)=0$, so $$f(0)=b_0+b_1(0)+b_2(0)^2+b_3(0)^3+b_4(0)^4+b_5(0)^5+...=0$$ and $b_0=0$.\\

Using the polynomial expression above, we obtain the following equations:\\

$$1-\frac{f(t)^2}{2}=(1-\frac{b_0^2}{2})-(b_0b_1)t-(b_0b_2+\frac{b_1^2}{2})t^2-(b_0b_3+b_1b_2)t^3-...$$

$$-\frac{a}{2}\biggl(1-\frac{f(t)^2}{2}\biggl)=-\frac{a}{2}\biggl((1-\frac{b_0^2}{2})-(b_0b_1)t-(b_0b_2+\frac{b_1^2}{2})t^2-(b_0b_3+b_1b_2)t^3-...\biggl)$$

$$-\frac{a}{2}\biggl(1-\frac{f(t)^2}{2}\biggl)=\bigl(-\frac{a}{2}+\frac{a}{4}b_0^2\bigl)+\bigl(\frac{a}{2}b_0b_1\bigl)t+\bigl(\frac{a}{4}b_1^2+\frac{a}{2}b_0b_2\bigl)t^2+\bigl(\frac{a}{2}b_0b_3+\frac{a}{2}b_1b_2\bigl)t^3+...$$

Now, we can match the coefficients of $f''(t)$ and $-\frac{a}{2}\bigl(1-\frac{f(t)^2}{2}\bigl)$ and solve for them.

$$2b_2+6b_3t+12b_4t^2+...=\bigl(-\frac{a}{2}+\frac{a}{4}b_0^2\bigl)+\bigl(\frac{a}{2}b_0b_1\bigl)t+\bigl(\frac{a}{4}b_1^2+\frac{a}{2}b_0b_2\bigl)t^2+...$$

$$-\frac{a}{2}+\frac{a}{4}b_0^2=2b_2$$

$$\frac{a}{2}b_0b_1=6b_3$$

$$\frac{a}{4}b_1^2+\frac{a}{2}b_0b_2=12b_4$$

$$\vdots$$

Once we know the values of $b_0$ and $b_1$, we can find the values for all of the rest of the coefficients $b_n$. Although we already know that $b_0=0$, we can't find $b_1$ from the system of equations above, so we will just have to define $b_1$ as $k$, and solve for all of the other coefficients symbolically, in terms of $k$ and $a$.\\

We obtain the following values for the coefficients:
$$b_0=0$$
$$b_1=k$$
$$b_2=-\frac{a}{4}$$
$$b_3=0$$
$$b_4=\frac{a}{48}k^2$$
$$b_5=-\frac{a^2}{160}k$$
$$b_6=\frac{a^3}{1920}$$
$$b_7=\frac{a^2}{4032}k^3$$
$$b_8=-\frac{11a^3}{107520}k^2$$
$$b_9=\frac{a^4}{69120}k$$
$$b_{10}=\frac{a^3}{387072}k^4-\frac{a^5}{1382400}$$
$$\vdots$$

Therefore, the polynomial becomes $$f(t)=(k)t-(\frac{a}{4})t^2+(\frac{a}{48}k^2)t^4-(\frac{a^2}{160}k)t^5+(\frac{a^3}{1920})t^6+(\frac{a^2}{4032}k^3)t^7+...$$

Using the condition that $f'(1)=0$ (i.e., the sum of the coefficients of $f(t)$ is 0), we found that a reasonable value for $b_1=k$ when $a=0.5$ was $k=0.249$. We then graphed the polynomial approximation of $f(t)$ in Desmos:

\begin{figure}[H]
    \centering
    \includegraphics[width=0.9\linewidth]{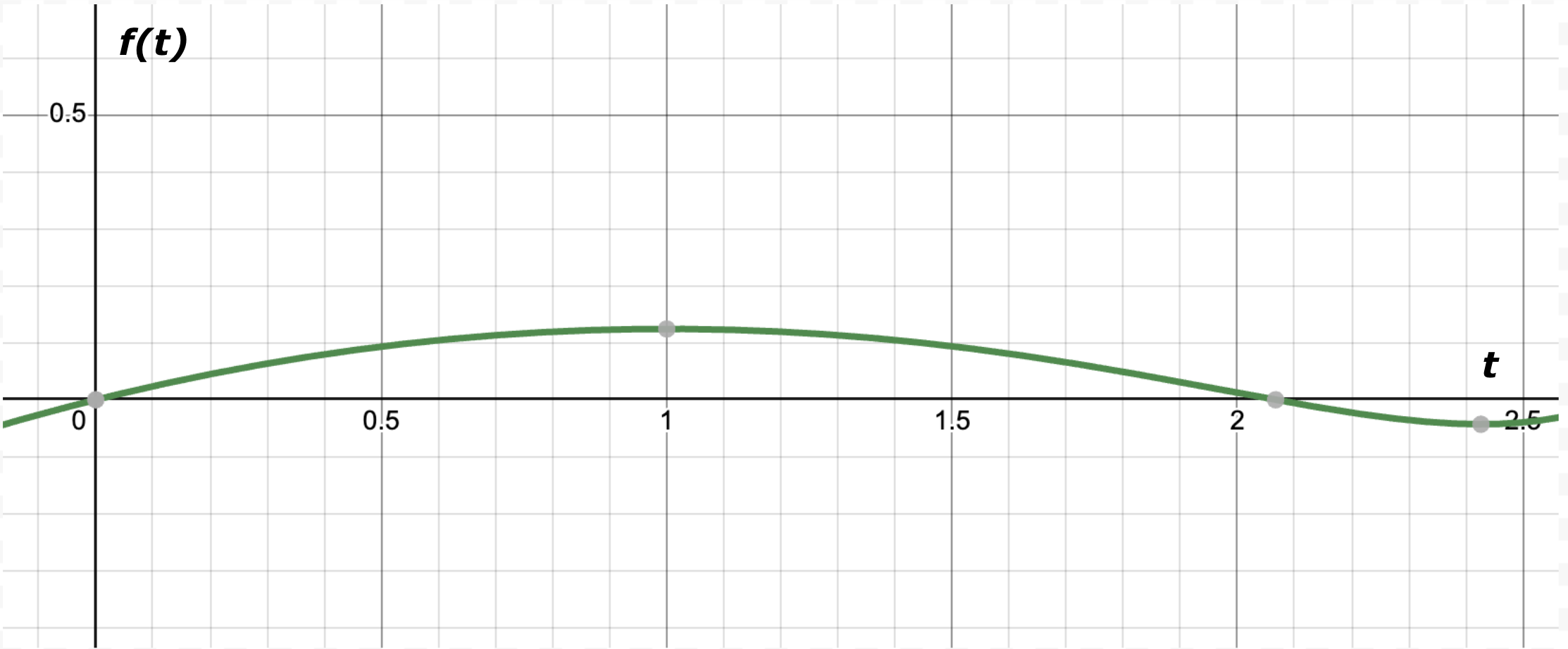}
    \caption{Polynomial approximation of $f(t)$ with 10 terms for when $a=0.5$ and $k=0.249.$}
    \label{fig:polyapprox}
\end{figure}

The function seems to fit our initial conditions of $f(0)=f'(1)=0$. However, if we increase $a$ to around $a=5$, and then use Desmos and Excel to find a reasonable value for $k$, the polynomial approximation is not as accurate. This is revealed by the increasing percent difference values in Table 1, which we obtained after comparing the $f''(1)$ approximation values from our polynomial function with the actual $f''(1)$ values, which are $-\frac{a}{2}\cos(f(1))$, from our differential equation.\\

\begin{table}[H]
    \centering
    \begin{tabular}{ |c|c|c|c|c|} 
\hline
$a$ & $k$ & $f''(1)$ approximation & $f''(1)=-\frac{a}{2}\cos(f(1))$ & \% difference\\
\hline
0.5 & 0.249 & -0.2481 & -0.2481 & 0.001\\
\hline
1 & 0.492 & -0.4850 & -0.4851 & 0.018\\
\hline
2 & 0.943 & -0.8925 & 0.8955 & 0.340\\
\hline
3 & 1.334 & -1.181 & -1.207 & 2.168\\
\hline
4 & 1.665 & -1.314 & -1.437 & 8.551\\
\hline
5 & 1.942 & -1.221 & -1.625 & 24.835\\
\hline
\end{tabular}
    \caption{Accuracy of polynomial approximation for different $a$ values. $k$ is the estimated value that satisfies the polynomial approximation above. The third column is the value from the Taylor polynomial approximation while the fourth column is the value from an Euler's approximation with a step size of $dt = 0.000001$. The approximation is meant for values of $a<<1$ and the table demonstrates why our approximation fails for larger values of $a$.}
    \label{tab:my_label}
\end{table}

Since we only used two terms in our Taylor approximation for $\cos(f(t))$, and we are multiplying the Taylor approximation by $a$, a larger $a$ value will result in a larger error. It seems that the strategy for minimizing $S$ depends on what the $a$ values are. We hypothesize that there are different Nash Equilibria for small, intermediate, and large $a$. We go into more detail about the Nash equilibrium for small $a$ in section 6.

%% file: 5_TwoLemmas.tex
\section{Two Lemmas}
\subsection{Comparing $y$ and $\sin(y)$}
\textbf{Conjecture:}
$$ky^2\ge y-\sin(y) \text{ for any } k \in [\frac{1}{\pi}, \infty).$$

\textbf{Proof:} Let us consider the expression $$p(y)=\frac{y-\sin(y)}{y^2}=\frac{1}{y}-\frac{\sin(y)}{y^2}.$$

To find where $p(y)$ has a maximum, we set the derivative of $p(y)$ equal to 0.

$$p'(y)=-\frac{1}{y^2}-\bigg(\frac{y^2\cos(y)-\sin(y)(2y)}{y^4}\bigg)=\frac{2\sin(y)-y\cos(y)-y}{y^3}=0$$

This occurs when $y=(2k+1)\pi$ and $k$ is an integer, which is where $p(y)$ has critical points. However, we want to find when $p(y)$ has an absolute maximum, so we will find the second derivative and find when it is negative.\\

The second derivative of $p(y)$ is as follows:
$$p''(y)=\frac{2y+4y\cos(y)-6\sin(y)+y^2\sin(y)}{y^4}.$$

For all $y=(2k+1)\pi$, $\cos(y)=-1$ and $\sin(y)=0$. Therefore, the second derivative becomes
$$p''(y)=\frac{2y-4y}{y^4}=\frac{-2}{y^3}.$$

For all $k\ge0$, $y=(2k+1)\pi>0$, and $p''(y)=\frac{-2}{y^3}<0$. This means that for positive $k$ values, $p(y)$ will have a local maximum.\\

For all $k<0$, $y=(2k+1)\pi<0$, and $p''(y)=\frac{-2}{y^3}>0$. This means that for all negative $k$ values, $p(y)$ will have a local minimum.\\

Now, we will focus on all $y=(2k+1)\pi>0$ with $k\ge0$.
$$p(y)=\frac{1}{y}-\frac{\sin(y)}{y^2}=\frac{1}{y}$$

$$p(\pi)=\frac{1}{\pi}>p(3\pi)=\frac{1}{3\pi}>p(5\pi)=\frac{1}{5\pi}>\cdots$$

Therefore, $p(y)$ has an absolute maximum when $y=\pi$, and the maximum value of $p(y)$ is $\frac{1}{\pi}\approx0.318$.

\begin{figure}[H]
    \centering
    \includegraphics[width=1\linewidth]{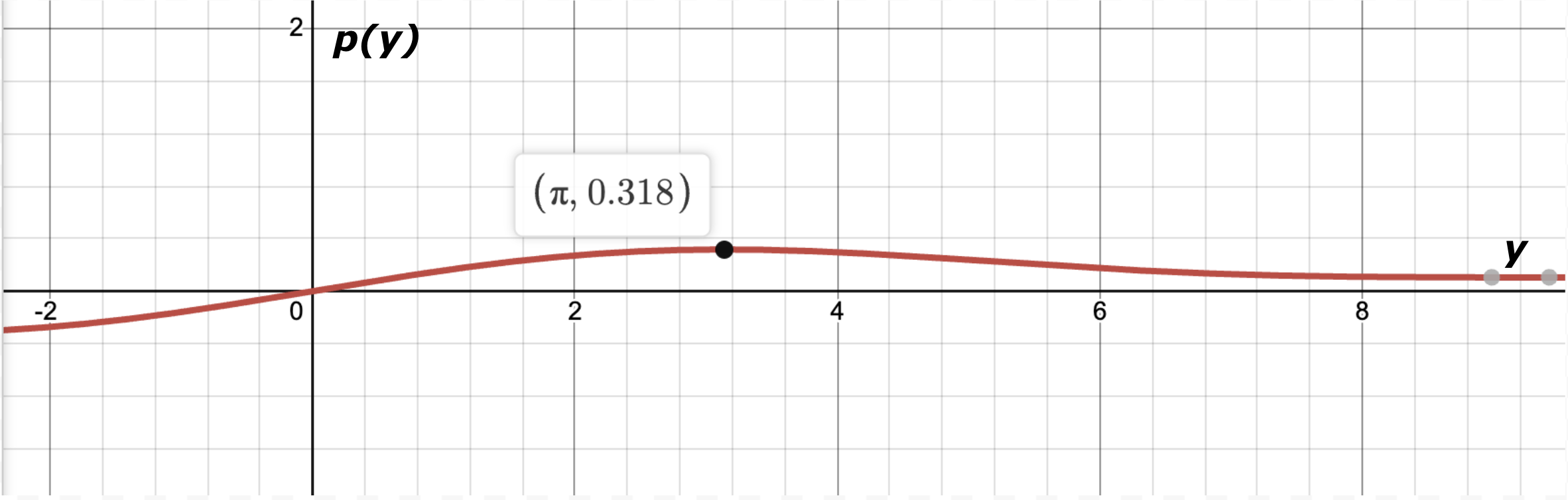}
    \caption{Graph of $p(y)$, showing absolute maximum at $y=\pi$.}
    \label{fig:p(y)max}
\end{figure}

We state that $$p(y)=\frac{y-\sin(y)}{y^2}\le\frac{1}{\pi}.$$

In other words, $$p(y)=\frac{y-\sin(y)}{y^2}\le k \qquad \forall \ k\ge\frac{1}{\pi}.$$

Rearranging, we obtain the expression $$ky^{2}\ge y-\sin(y)  \qquad \forall \ k\ge\frac{1}{\pi}.$$

\subsection{Integrals and Fourier Series}
\textbf{Conjecture:}
$$\int_{0}^{2}f'(t)^2\ dt\ge\bigg(\frac{\pi}{2}\bigg)^2\int_{0}^{2}f(t)^2\ dt$$
for any differentiable function $f(t)$ where $f(0)=f(2)=0$.\\\\

\textbf{Proof:}
Let $f(t)$ be any function such that $f(0)=0$ and $f(2)=0$. We are only concerned with the function $f(t)$ on the interval $x \in [0,2]$, and so we can manipulate $f(t)$ such that it is an odd function. This determines how $f(t)$ behaves on the interval $[-2,2]$. Next, we will make $f(t)$ periodic with period 4 by repeating itself. Therefore, we can state that $f(t)$ is an odd function with period 4. Notice that this claim is true regardless of the shape of $f(t)$ on the interval $[0,2]$.\\

Now, we can express $f(t)$ with a Fourier Series since it is periodic. For any odd function with period $2L$, the Fourier Series expression is as follows:
$$f(t)=a_1\sin\bigg(\frac{\pi}{L}t\bigg)+a_2\sin
\bigg(\frac{2\pi}{L}t\bigg)+a_3\sin\bigg(\frac{3\pi}{L}t\bigg)+...$$
$f(t)$ has a period of 4, so $2L=4$ and $L=2$. Thus,
$$f(t)=a_1\sin\bigg(\frac{\pi}{2}t\bigg)+a_2\sin(\pi t)+a_3\sin\bigg(\frac{3\pi}{2}t\bigg)+...$$

Then, we can write the following equations:
$$f(t)^2=a_1^2\sin^2\bigg(\frac{\pi}{2}t\bigg)+2a_1a_2\sin\bigg(\frac{\pi}{2}t\bigg)\sin(\pi t)+...$$

$$f'(t)=\bigg(\frac{\pi}{2}\bigg)a_1\cos\bigg(\frac{\pi}{2}t\bigg)+\big(\pi\big) a_2\cos(\pi t)+\bigg(\frac{3\pi}{2}\bigg)a_3\cos\bigg(\frac{3\pi}{2}t\bigg)+...$$

$$f'(t)^2=\bigg(\frac{\pi}{2}\bigg)^2a_1^2\cos^2\bigg(\frac{\pi}{2}t\bigg)+2\bigg(\frac{\pi}{2}\bigg)\big(\pi\big)a_1a_2\cos\bigg(\frac{\pi}{2}t\bigg)\cos(\pi t)+...$$\\

Now, we will take the integral of $f(t)^2$ and $f'(t)^2$ on the interval $[0,2]$.\\

Due to the orthogonal relationships of sine and cosine functions, any terms in the integral of the form $\sin(mx)\sin(nx)$ or $\cos(mx)\cos(nx)$, with $m\neq n$, are equal to zero \cite{fouriersin}. So,

$$\int_{0}^2f(t)^2\ dt=\int_{0}^2\bigg[a_1^2\sin^2\bigg(\frac{\pi}{2}t\bigg)+a_2^2\sin^2(\pi t)+a_3^2\sin^2\bigg(\frac{3\pi}{2}t\bigg)+...\bigg]\ dt$$

$$\int_{0}^2f'(t)^2\ dt=\int_{0}^2\bigg[\bigg(\frac{\pi}{2}\bigg)^2a_1^2\cos^2\bigg(\frac{\pi}{2}t\bigg)+(\pi)^2a_2^2\cos^2(\pi t)+\bigg(\frac{3\pi}{2}\bigg)^2a_3^2\cos^2\bigg(\frac{3\pi}{2}t\bigg)+...\bigg]\ dt.$$\\

Note that in $f(t)^2$, all of the squared sines are in the form $\sin^2\big(\frac{\pi}{2}kt\big)$, where $k$ is a positive integer. For any term $\sin^2\big(\frac{\pi}{2}kt\big)$, $$\int_{0}^{2}\sin^2\big(\frac{\pi}{2}kt\big)\ dt=\int_{0}^{2}\frac{1-\cos(k\pi t)}{2}\ dt=1, \ \forall k.$$

The same can be said for the cosine expressions in $f'(t)^2$, since all of them are in the form $\cos^2\big(\frac{\pi}{2}kt\big)$. For any of these terms, $$\int_{0}^{2}\cos^2\big(\frac{\pi}{2}kt\big)
\ dt=\int_{0}^{2}\frac{1+\cos(k\pi t)}{2}\ dt=1, \ \forall k.$$

Plugging in these observations, we get $$\int_{0}^{2}f(t)^2\ dt=a_{1}^{2}+a_{2}^{2}+a_{3}^{3}+...$$

$$\int_{0}^{2}f'(x)^2\ dt=\bigg(\frac{\pi}{2}\bigg)^2a_1^2+(\pi)^2a_2^2+\bigg(\frac{3\pi}{2}\bigg)^2a_3^2+...$$\\

Now, notice that we can make the following claim:
$$\bigg(\frac{\pi}{2}\bigg)^2a_1^2+\big(\pi\big)^2a_2^2+\bigg(\frac{3\pi}{2}\bigg)^2a_3^2+...\ge\bigg(\frac{\pi}{2}
\bigg)^2a_{1}^{2}+\bigg(\frac{\pi}{2}
\bigg)^2a_{2}^{2}+\bigg(\frac{\pi}{2}
\bigg)^2a_{3}^{3}+...$$

$$\int_{0}^{2}f'(x)^2\ dt\ge\bigg(\frac{\pi}{2}
\bigg)^2\int_{0}^{2}f(t)^2\ dt.$$

The expression above only applies to odd functions $f(t)$ with period 4 where $f(0)=f(2)=0$.

%% file: 6_SmallA.tex
\section{When $a$ is small}
Since the differential equation $$f''(t)=-\frac{a}{2}\cos\bigl(f(t)-g(t)\bigl)$$ is coupled with the opponent's move $g(t)$, there is no easy Nash equilibrium, or optimal strategy, for all values of $a$. So, for this section, we will focus on the case where $a$ is small first.

\subsection{Calculus of Variations for Small Values of $a$}
Recall that the Taylor polynomial for $\sin(x)$ starts with $x$. When $a$ is very small, we can make the following approximation: $a\sin\bigl(f(t)-g(t)\bigl)\  \approx a\big(f(t)-g(t)\big)$, by assuming that the rest of the higher order terms of the Taylor series will  become negligible after being multiplied by a very small $a$ value.\\

We can rewrite $S$ as $$S=\int_{0}^{1}\bigg(f'(t)^2-g'(t)^2-a\sin\bigl(f(t)-g(t)\bigl)\bigg)\ dt\approx\int_{0}^{1}\bigg(f'(t)^{2}-g'(t)^{2}-a\bigl(f(t)-g(t)\bigl)\bigg)\ dt.$$\\

Now, we will apply calculus of variations again to obtain more information on what the optimal function is. \cite{calcofvariations}

$$0=\diff{}{s}S\bigg\rvert_{s=0}=\dot{S}\bigg\rvert_{s=0}$$
$$=\int_{0}^{1}\frac{\partial}{\partial s}\biggl(f'(t)^{2}-g'(t)^{2}-a\bigl(f(t)-g(t)\bigl)dt\biggl)\bigg\rvert_{s=0}dt$$
$$=\int_{0}^{1}\biggl(2f'(t)\dot{f}'(t)-a\dot{f}(t)\biggl)\bigg\rvert_{s=0}dt.$$\\

Again, using integration by parts, and the fact that $f'(1)=0$ and $\dot{f}(0)=0$, we find that
$$\int_{0}^{1}f'(t)\dot{f}'(t)\ dt =\int_{0}^{1}-f''(t)\dot{f}(t)\ dt.$$\\

Plugging this in, we get that
$$0=\int_{0}^{1}\bigg(-2f''(t)\dot{f}(t)-a\dot{f}(t)\bigg)\bigg\rvert_{s=0}dt=\int_{0}^{1}\bigl(-2f''(t)-a\bigl)\dot{f}(t)\bigg\rvert_{s=0}dt $$

$$=\int_{0}^{1}\bigl(-2F''(t)-a\bigl)\eta(t)\ dt \qquad \quad \forall \eta(t).$$

This implies that $$-2F''(t)-a=0$$

$$F''(t)=-\frac{a}{2}$$

$$F'(t)=-\frac{a}{2}t+b.$$

Now, we can use our initial condition that $F'(1)=0.$

$$F'(1)=-\frac{a}{2}+b=0 \Rightarrow b=\frac{a}{2}$$

$$F'(t)=-\frac{a}{2}t+\frac{a}{2}$$

$$F(t)=-\frac{a}{4}t^{2}+\frac{a}{2}t+c$$\\

Plugging in the initial condition that $F(0)=0$, we get that $F(0)=c=0.$\\

So, our final equation is $$F(t)=-\frac{a}{4}t^{2}+\frac{a}{2}t=-\frac{a}{4}(t^2-2t)=-\frac{a}{4}t(t-2)=\frac{a}{4}t(2-t).$$

To satisfy the condition that $F''(t)=-\frac{a}{2}$ and the condition that $F(0)=F'(1)=0$, there is only 1 function that exists, which is shown above.\\

Now, we must prove that this function is the optimal function. We know that when $F(t)=\frac{a}{4}t(2-t)$, the functional $S$ has a critical point, but we haven't proven that this critical point is a minimum yet.\\

To prove that $F(t)=\frac{a}{4}t(2-t)$ is a minimum of the functional $S$ at $s=0$, we evaluate the second derivative of $S$ at $s=0$.

$$\diff{^2}{s^2}S\bigg\rvert_{s=0}=\diff{}{s}\dot S\bigg\rvert_{s=0}=\int_{0}^{1}\frac{\partial}{\partial s}\biggl(2f'(t)\dot{f}'(t)-a\dot{f}(t)\biggl)\bigg\rvert_{s=0}dt$$

$$=\int_{0}^{1}\frac{\partial}{\partial s}\biggl(2f'(t)\eta'(t)-a\eta(t)\biggl)\bigg\rvert_{s=0}dt=\int_{0}^{1}2\eta'(t)^2\bigg\rvert_{s=0}dt=\int_{0}^{1}2\eta'(t)^2dt>0$$\\

The second derivative of $S$ at $s=0$ is always greater than 0, which means that $F(t)=\frac{a}{4}t(2-t)$ is a minimum. Since there was only one critical point, $F(t)=\frac{a}{4}t(2-t)$ must also be the \textbf{absolute} minimum.

$$f_{optimal}(t)=F(t)=\frac{a}{4}t(2-t)$$

$$S(f_{optimal}(t))=-\frac{a^{2}}{12} \text{ when } g(t)=0$$

\subsection{Proof of Nash equilibrium for Small Values of $a$}

\textbf{Conjecture: }For sufficiently small values of $a$, there is only one optimal solution for the Nash equilibrium. This optimal solution is
$$f(t)=\frac{a}{4}t(2-t), \ g(t)=\frac{a}{4}t(2-t).$$

We have already proven that the optimal function for player $f$ to play is $f(t)=\frac{a}{4}t(2-t)$, so now, we have to prove that the most optimal function for $g$ to play is also $g(t)=\frac{a}{4}t(2-t)$.\\

Since $S\big(f(t),g(t)\big)$ is an odd functional, when $f(t)=g(t)=\frac{a}{4}t(2-t)$, the value of the functional will be zero: $S\big(f(t),f(t)\big)=-S\big(f(t),f(t)\big)=0$.\\

Recall that player $g$ wants to maximize $S\big(f(t),g(t)\big)$, or minimize $S\big(g(t),f(t)\big)$. We want to show that for any  function that player $g$ chooses, which we will denote by $g(t)$, $$S\big(g(t),f(t)\big)\ge S\big(f(t),f(t)\big)=0.$$

In other words, the most optimal function for player $g$ to play, or the solution that will minimize $S\big(g(t),f(t)\big)$, is $g(t)=f(t)=\frac{a}{4}t(2-t)$.\\\\

\textbf{Proof: }Let $g(t)$ be defined as $f(t)$ with some slight variation, which we will denote as $h(t)$. So, $g(t)=f(t)-h(t)$. We then can substitute $f(t)-h(t)$ in for $g(t)$ in the functional.

$$S\big(g(t),f(t)\big)=S\big(f(t)-h(t),f(t)\big)$$\\

We want to prove that
$$S\big(f(t)-h(t),f(t)\big) \ge S\big(f(t),f(t)\big)=0.$$

Next, we expand the functional into its full form:
$$S\big(f(t)-h(t),f(t)\big)=\int_{0}^{1}\bigg(\big(f(t)-h(t)\big)'^2-f'(t)^2-a\sin\bigl(\big(f(t)-h(t)\big)-f(t)\bigl)\bigg)\ dt$$

$$=\int_{0}^{1}\bigg(\big(f'(t)-h'(t)\big)^2-f'(t)^2-a\sin(-h(t))\bigg)\ dt$$

$$=\int_{0}^{1}\bigg(\big(f'(t)^2-2f'(t)h'(t)+h'(t)^2\big)-f'(t)^2-a\sin(-h(t))\bigg)\ dt$$

$$=\int_{0}^{1}\bigg(-2f'(t)h'(t)+h'(t)^2+a\sin(h(t))\bigg)\ dt.$$\\

Now, we will use integration by parts to simplify the first term of the integral.

$$\int_{0}^{1}f'(t)h'(t)\ dt=f'(t)h(t)\bigg|_{0}^{1}-\int_{0}^{1}f''(t)h(t)\ dt$$

However, recall that $f(0)=f'(1)=0$ and $h(0)=h'(1)=0$, so 

$$\int_{0}^{1}f'(t)h'(t)\ dt=-\int_{0}^{1}f''(t)h(t)\ dt.$$\\

Plugging this in, we get
$$S\big(f(t)-h(t),f(t)\big)=\int_{0}^{1}\bigg(2f''(t)h(t)+h'(t)^2+a\sin(h(t))\bigg)\ dt.$$\\

Using the fact that $f''(t)=-\frac{a}{2}$ since $f(t)=\frac{a}{4}t(2-t)$, we get

$$S\big(f(t)-h(t),f(t)\big)=\int_{0}^{1}\bigg(-ah(t)+h'(t)^2+a\sin(h(t))\bigg)\ dt.$$\\

Now, we will use the relationship between $\sin(h)$ and $h$ that we proved in Section 5.1. From previous result, we know that
$$kh^{2}\ge h-\sin(h) \qquad \forall \ k\ge\frac{1}{\pi}.$$

So, the following must be true: $$\frac{1}{\pi}h^{2}\ge h-\sin(h)$$
$$\sin(h)\ge h-\frac{1}{\pi}h^{2}.$$

Subbing in $\sin(h(t))$, we get

$$\int_{0}^{1}\bigg(-ah(t)+h'(t)^2+a\sin(h(t))\bigg)\ dt \ge \int_{0}^{1}\bigg(-ah(t)+h'(t)^2+a\big[h(t)-\frac{1}{\pi}h(t)^{2}\big]\bigg)\ dt.$$

Notice that the $ah(t)$ terms cancel out nicely and we are left with

$$\int_{0}^{1}\bigg(-ah(t)+h'(t)^2+a\sin(h(t))\bigg)\ dt \ge \int_{0}^{1}\bigg(h'(t)^2-\frac{a}{\pi}h(t)^{2}\bigg)\ dt.$$

The expression on the right is similar to the conjecture that we proved using integrals and Fourier Series in Section 5.2. However, before we use it, we need to check to see if $h(t)$ meets the conditions.\\

In the game, the only initial conditions given are that $h(0)=h'(1)=0$. However, in order to use the conjecture from Section 5.2, we also need $h(2)=0$.\\

Let $h(t)$ be any function on the interval $t \in [0,1]$. Then, on the interval $t \in (1,2]$, define $h(t)$ such that it is symmetrical about the line $t=1$. In other words, it is a reflection of itself from $[0,1]$ about the line $t=1$. Since we are given that $h(0)=0$, $h(2)$ must also equal 0 (purely due to how we defined $h(t)$). So, we have now adjusted $h(t)$ such that it fits the required conditions, and we can use the conjecture from Section 5.2 on $h(t)$.\\

From Section 5.2, we know that $$\int_{0}^{2}h'(t)^2\ dt \ge \int_{0}^{2}\big(\frac{\pi}{2}\big)^2h(t)^2\ dt.$$\\

Since $h(t)$ is symmetric about the line $t=1$, $h(t)^2$ and $h'(t)^2$ are also symmetric about the line $t=1$. So, we can adjust the limits of integration for our statement. 
$$\int_{0}^{2}h'(t)^2\ dt =2\int_{0}^{1}h'(t)^2\ dt \ge \int_{0}^{2}\big(\frac{\pi}{2}\big)^2h(t)^2\ dt=2\int_{0}^{1}\big(\frac{\pi}{2}\big)^2h(t)^2\ dt$$

$$\int_{0}^{1}h'(t)^2\ dt \ge \int_{0}^{1}\big(\frac{\pi}{2}\big)^2h(t)^2\ dt.$$\\

Now, we can substitute the inequality into our original statement.

$$\int_{0}^{1}\bigg(h'(t)^2-\frac{a}{\pi}h(t)^{2}\bigg)\ dt \ge \int_{0}^{1}\bigg(\big(\frac{\pi}{2}\big)^2h(t)^2-\frac{a}{\pi}h(t)^{2}\bigg)\ dt=\int_{0}^{1}\bigg(\bigg[\big(\frac{\pi}{2}\big)^2-\frac{a}{\pi}\bigg]h(t)^2\bigg)\ dt$$

Remember we wanted to show that $S\big(f(t)-h(t),f(t)\big)\ge0$. This is only true for certain small values of $a$.

$$S\big(f(t)-h(t),f(t)\big)\ge\int_{0}^{1}\bigg(\bigg[\big(\frac{\pi}{2}\big)^2-\frac{a}{\pi}\bigg]h(t)^2\bigg)\ dt$$

In order for the expression on the right-hand side to be greater than or equal to 0, $\big(\frac{\pi}{2}\big)^2-\frac{a}{\pi}$ must be greater than or equal to 0.

$$\big(\frac{\pi}{2}\big)^2-\frac{a}{\pi}\ge0$$

Solving for $a$, we obtain the following inequality:

$$\frac{a}{\pi}\le\big(\frac{\pi}{2}\big)^2=\frac{\pi^2}{4}$$

$$a\le\frac{\pi^3}{4}.$$

Since $a$ is a nonnegative number, $$0\le a\le\frac{\pi^3}{4}\approx7.75.$$

Notice that if we had chosen a larger $k$ value in our comparison of $\sin(h(t))$ and $h^2$, then we wouldn't obtain the upper bound on $a$ when $f(t)=g(t)=\frac{a}{4}t(2-t)$ is the Nash equilibrium. However, we wanted an upper bound on how small $a$ must be, so using $k=\frac{1}{\pi}$ gives the best result.\\

In conclusion, for all $a \in [0,\frac{\pi^3}{4}]$, the Nash equilibrium to the functional $S$ is $f(t)=g(t)=\frac{a}{4}t(2-t)$.

%% file: 7_LargeA.tex
\section{When $a$ is large}
Next, we approached the problem of finding the Nash equilibrium when $a$ is large (i.e. for $a\ge \frac{\pi^3}{4}$). We hypothesized that the Nash equilibrium for larger values of $a$ will become a mixed strategy, specifically one with three optimal strategies, that are each played a certain percentage of the time, and operates like a weighted rock-paper-scissors game.

\subsection{Explanation of Linear Programming}
In our procedure for finding the Nash equilibrium for larger values of $a$, we implement a process called linear programming. Before we dive into what our procedure is, we will first explain the process of linear programming.\\

Take a = 75 as an example. Let us consider the 3 functions $r(t), p(t)$, and $s(t)$, which are shown below.\\

\begin{figure}[H]
    \centering
    \includegraphics[width=1\linewidth]{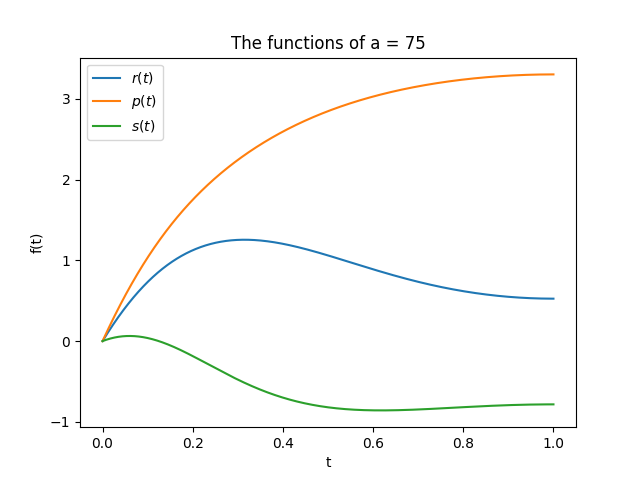}
    \caption{Graph of the example functions $p(t), s(t),$ and $r(t)$, for when $a=75$.}
    \label{fig:a=75}
\end{figure}

The goal of linear programming is to find out what the optimal strategy is if we are only allowed to play these three functions. In other words, how often should you play the functions $r(t), p(t)$, and $s(t)$ such that you maximize your chances of winning? In order to find these probabilities, we create a payoff matrix, where the payoff values in the matrix represent how much player $f$ will win if they play the function of that column against player $g$, who is playing the function of that row. The setup of this matrix is shown below:\\

\[
\begin{blockarray}{cccc}
& r(t) & p(t) & s(t)\\
\begin{block}{c(ccc)}
  r(t) & S(r,r) & S(p,r) & S(s,r)\\
  p(t) & S(r,p) & S(p,p) & S(s,p)\\
  s(t) & S(r,s) & S(p,s) & S(s,s)\\
\end{block}
\end{blockarray}
 \]

Once we plug in the functions $r(t), p(t), $ and $s(t)$ into our payoff matrix, it looks like this:

\[
\begin{blockarray}{cccc}
& r(t) & p(t) & s(t)\\
\begin{block}{c(ccc)}
  r(t) & 0 & 34.994 & -62.147\\
  p(t) & -34.994 & 0 & 32.740\\
  s(t) & 62.147 & -32.740 & 0\\
\end{block}
\end{blockarray}
 \]
 
Next, we will define $P_r$, $P_p$, and $P_s$ to be the probability that Player $f$ plays the functions $r(t)$, $p(t)$ and $s(t)$ for the Nash equilibrium strategy, respectively. We can state the following about these probabilities:
$$P_r+P_p+P_s=1$$
$$P_r,P_p,P_s\ge0.$$

Now, assume that Player $g$ is playing $r(t)$. Then, the expected value for the amount of money Player $f$ would earn is $$0\ge E_r=(0)\cdot P_r+(34.994)\cdot P_p+(-62.147)\cdot P_s.$$ Notice that this expected value must be less than or equal to 0 (Assume you are player 1 and player 2 is playing the Nash equilibrium; then, you want to maximize your payoff by having the value be less than or equal to 0). Similarly, we can define $E_p$ and $E_s$, which correspond to when Player $g$ plays $p(t)$ and $s(t)$, respectively. $$0\ge E_p=(-34.994)\cdot P_r+(0)\cdot P_p+(32.740)\cdot P_s$$ $$0\ge E_s=(62.147)\cdot P_r+(-32.740)\cdot P_p+(0)\cdot P_s$$

Using matrices to simplify this system of linear equations, we get\\
\[ \begin{bmatrix}
         0\\
         0\\ 
         0\\ 
     \end{bmatrix}
     \ge
     \begin{bmatrix}
         0 & 34.994 & -62.147\\
         -34.994 & 0 & 32.740\\ 
         62.147 & -32.740 & 0\\
     \end{bmatrix}
     \begin{bmatrix}
         P_r\\
         P_p\\ 
         P_s\\
     \end{bmatrix} \]

From here, we can observe that the system of inequalities forms a cycle, since we have $$62.147P_s\ge34.994P_p$$
$$34.994P_r\ge32.740P_s$$
$$32.740P_p\ge62.147P_r.$$

Therefore, the statements only hold true if all the values are equal. So, we are able to solve for the probabilities $P_r$, $P_p$, and $P_s$ through some simple algebraic manipulation. For the example functions $r(t), p(t),$ and $s(t)$ above, the Nash equilibrium strategy is $(P_r, P_p, P_s)=(0.252 ,0.478 , 0.269)$.

\subsection{Procedure}

Now that we have explained the process of linear programming, we can dive into our procedure for finding the Nash equilibrium.

\begin{enumerate}
    \item Start with 6 seed functions that satisfy the initial conditions $f(0)=0$ and $f'(1)=0$. Call them $f_1(t), f_2(t), ..., f_6(t)$. Technically, we could have chosen to start with any number of functions, but the more functions we start with, the better, since they are more representative of all possible functions that $f(t)$ could be.
    \item Create a payoff matrix containing the values $S(f_i(t),f_j(t))$ for all $i,j \in \{1,2,3,...,6\}$, as shown previously.
    \item Run the payoff matrix through the linear programming code to find the probabilities that player $f$ plays each function $f_1(t), f_2(t), ..., f_6(t)$ for the Nash equilibrium strategy.
    \item Now, suppose that player $g$ is playing this Nash equilibrium strategy against player $f$, which we will denote as $(\{P_k\}^n_{k=1},\{g_k(t)\}^n_{k=1})$, where $P_k$ is the probability that player $g$ plays the function $g_k(t)$. We want to determine if there is a new optimal function player $f$ can play that will beat this Nash equilibrium strategy. We can find this new optimal function by rewriting $g(t)$ as a weighted sum of all of the functions in the Nash equilibrium strategy, where the weights are the probabilities. So, our expression for $S=S(f(t),g(t))$ becomes $$S=\sum^n_{k=1}\bigg[P_k\int_0^1 f'(t)^2-g_k'(t)^2-a\sin(f(t)-g_k(t))\ dt\bigg].$$ Using calculus of variations, we are able to minimize $S$ and solve for $f(t)$.
    \item Use Euler's method to determine if this new function beats the old strategy by more than the threshold of $10^{-5}$ (we will explain this threshold more in the next section).
    \item If this new function beats our old Nash equilibrium strategy by more than the set threshold, then we add it to our list of seed functions and repeat the process again (starting from Step $\#2$).
    \item If this new function does not beat our old Nash equilibrium strategy by more than $10^{-5}$, then we conclude that the old strategy is the Nash equilibrium and we are done.
\end{enumerate}

\begin{figure}[H]
    \centering
    \includegraphics[width=1\linewidth]{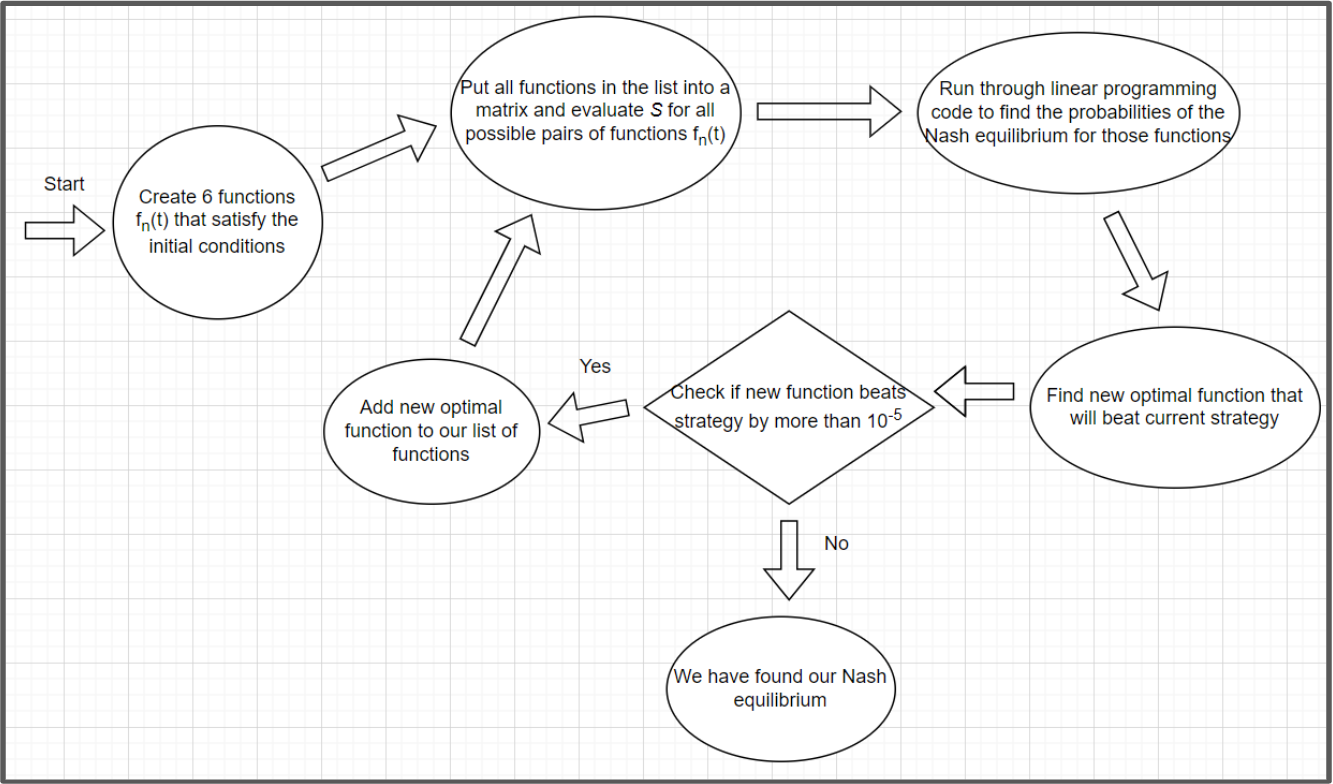}
    \caption{Flowchart outlining our process to find the Nash equilibrium.}
    \label{fig:flowchart}
\end{figure}

\subsection{Euler's Method, Riemann Sums}

Here’s a simplified explanation of the methodology for solving a complex differential equation without a fully defined point.\\

First off, we create an array of $f_{i}(t)$ values for each of the functions in our list (from $t=0$ to $t=1$). Since we have a second-order differential equation, we can  use Euler's method to solve for $f_{i}(0)$, since we want the function to satisfy the initial condition that $f_i(0)=0$ for all $i$.\\

Next, we look for the point where the $f(t)$ values at $t=0$ cross over the $x$-axis. This crossover point is crucial because it helps us identify a range of corresponding $f(t)$ values when $t=1$. Within this range, we use a search method to get $f(t)$ as close to $0$ as possible when $t=0$, because of the initial condition.\\

\href{https://www.desmos.com/calculator/npyu2ou7cv}{Here} is a link to a Desmos graph that illustrates how our method of solving a differential equation works.\\

Now that we have found all solutions to the differential equation that satisfy both initial conditions, we need to determine which one is the most optimal, since we are looking for the solution that produces the largest payoff. To do this, we test every possible function by plugging it back into the functional and comparing the payoff values.\\

To calculate the value of the functional $S(f(t), g(t)),$ we use a Riemann Sum, since our functions are defined as a list of points that we need to integrate over.\\

However, using these methods often lead to a notable lack of precision, where a mere difference of $10^{-5}$ can distinguish between one of the Nash equilibrium functions and its close estimates. Consequently, multiple approximate estimates of a Nash equilibrium function may cluster together. In such cases, we posit the existence of a singular true function amidst the array of estimated functions. Recognizing this inherent lack of precision, we experimented with smaller step sizes, which did aid in refining the solution. However, it's important to acknowledge that the solution remains an estimate regardless of our efforts.

\subsection{Examples for specific $a$ values}
After completing our code for our procedure for finding the Nash equilibrium, we tested various $a$ values and noticed that there was an interesting pattern. Examining the solutions provided below for $a$ values of $a = 10, 25, 50, 100, 175$, and $250$, a discernible pattern emerges regarding their relative positions. For $a = 10$ and $a = 25$, the Nash equilibrium is comprised of only two functions, deviating from our previously proposed three-function solution. However, this doesn't stay true for long, and the Nash equilibrium transitions from two functions to three functions at $a = 45$. Although $a$ values below $45$ exhibit a Nash equilibrium with only two functions, the same principle outlined earlier still applies when calculating the probabilities.\\

\begin{figure}[H]
    \centering
    \includegraphics[width=0.9\linewidth]{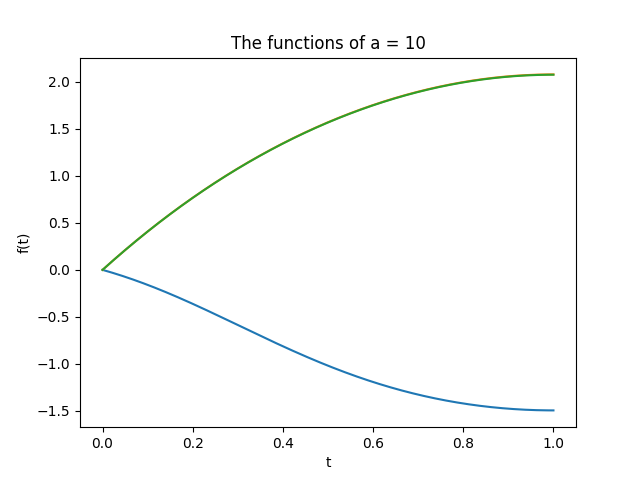}
    \caption{The two functions for the Nash equilibrium at $a = 10$.}
    \label{fig:a=10}
\end{figure}

\begin{figure}[H]
    \centering
    \includegraphics[width=0.9\linewidth]{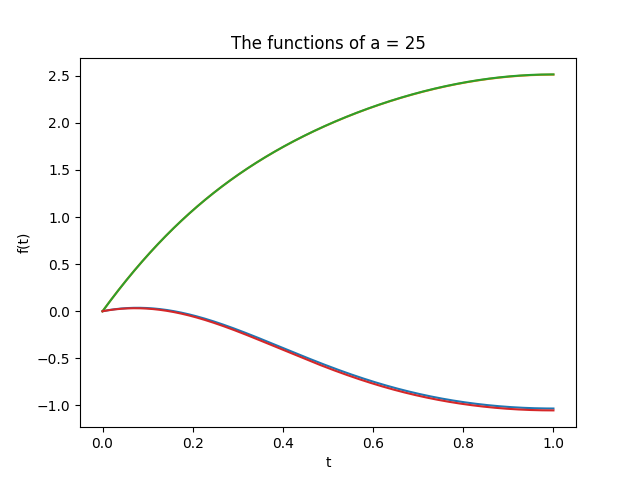}
    \caption{The two functions for the Nash equilibrium at $a = 25$.}
    \label{fig:a=25}
\end{figure}

\begin{figure}[H]
    \centering
    \includegraphics[width=0.9\linewidth]{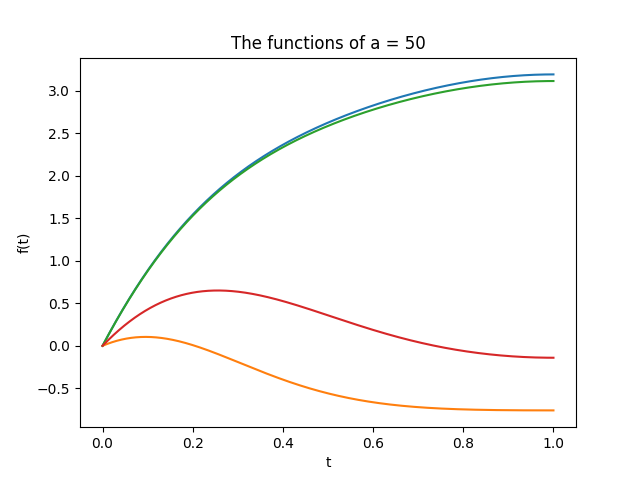}
    \caption{The three functions for the Nash equilibrium at $a = 50$.}
    \label{fig:a=50}
\end{figure}

\begin{figure}[H]
    \centering
    \includegraphics[width=0.9\linewidth]{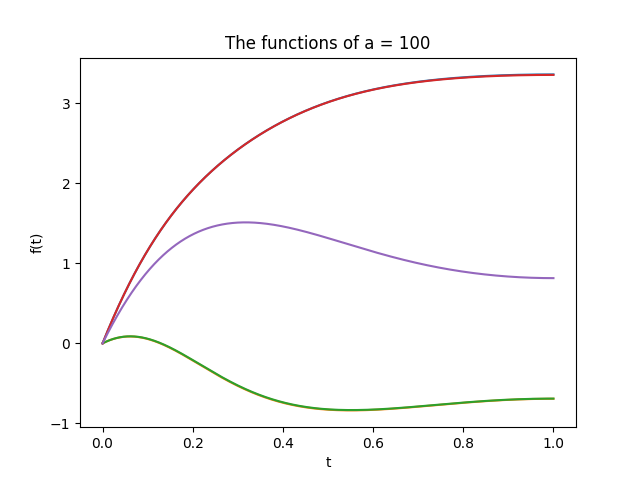}
    \caption{The three functions for the Nash equilibrium at $a = 100$.}
    \label{f100:a=100}
\end{figure}

\begin{figure}[H]
    \centering
    \includegraphics[width=0.9\linewidth]{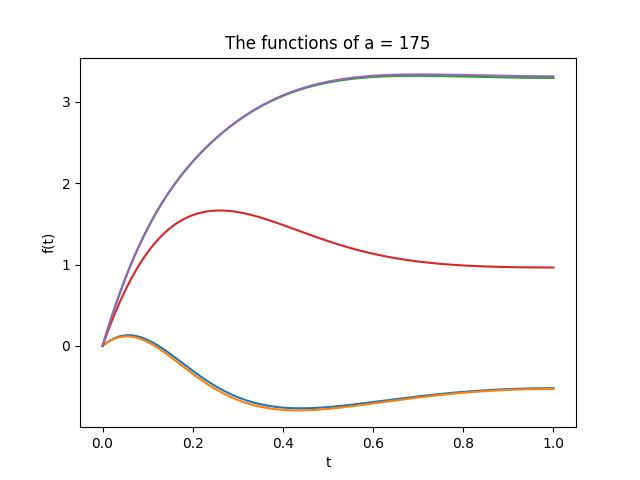}
    \caption{The three functions for the Nash equilibrium at $a = 175$.}
    \label{fig:a=175}
\end{figure}

\begin{figure}[H]
    \centering
    \includegraphics[width=0.9\linewidth]{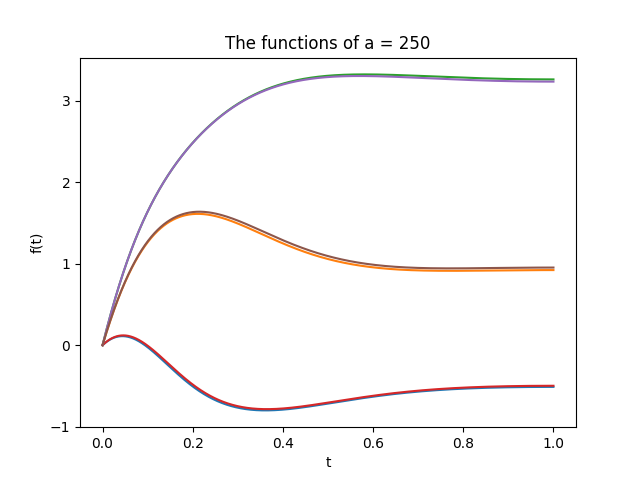}
    \caption{The three functions for the Nash equilibrium at $a = 250$.}
    \label{fig:a=250}
\end{figure}

Note: In the graphs for $a=25, 50, 175$, and $250$, you can see small clusters for a few of the Nash equilibrium functions. This occurs due to the lack of precision in Euler's Method and Riemann Sums.

\subsection{Results for $0\le a \le$ 256}
We became interested in the progression of the solution curves as the $a$ value increased. So, we decided to plot the endpoints of the solution curves against $a$, and observe when and where the Nash equilibrium changed from one function to two functions, and eventually three functions.
\begin{figure}[H]
    \centering
    \includegraphics[width=1\linewidth]{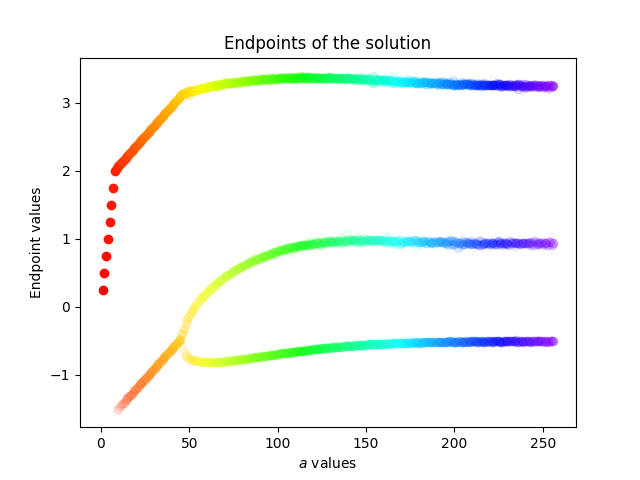}
    \caption{The progression of the endpoints of the functions that make up the Nash equilibrium solution from $a = 1$ to $a = 256$.}
    \label{fig:rainbowthingy}
\end{figure}

We were able to determine that the number of functions in the Nash equilibrium strategy jumps from two to three at $a=45$. This brought up the question of whether or not there would ever be three or more functions in the Nash equilibrium strategy.\\

Notice that in the graph above, the endpoints of the functions that make up the Nash equilibrium strategy have gradually stabilized as $a \rightarrow 256$. Therefore, we conjecture that the Nash equilibrium strategy will continue to be composed of three functions as $a$ increases. However, we have not proven this yet; this will be future work.

%% file: 8_sin_3x.tex
\section{Variations of the Game}
Recall that we could substitute the $\sin(t)$ function in the functional $S$ with any other odd function, like $\sin^3(t)$. So, a variation of the game is $$S=S(f(t),g(t))=\int_{0}^{1}\bigg(f'(t)^2-g'(t)^2-a\sin^3\big(f(t)-g(t)\big)\bigg)\ dt.$$ We decided to investigate what the progression of the Nash equilibrium solutions would look like if this was the case.

\begin{figure}[H]
    \centering
    \includegraphics[width=1\linewidth]{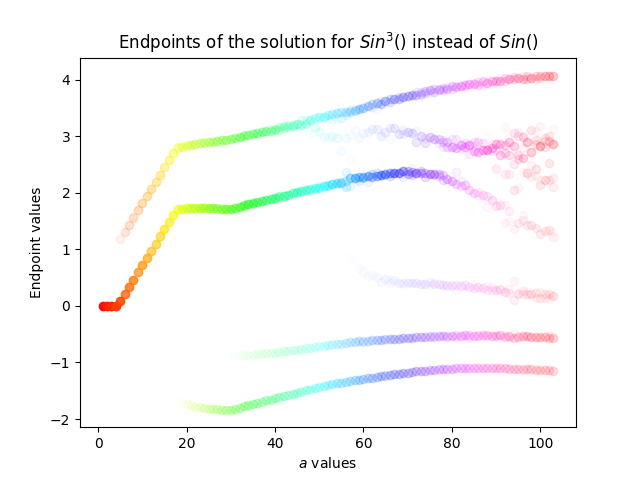}
    \caption{The progression of the endpoints of the functions that make up the Nash equilibrium solution from $a = 1$ to $a = 103$.}
    \label{fig:rainbowthingy2.0}
\end{figure}

In the $\sin^3(t)$ game, it initiates with a single solution similar to the $\sin(t)$ game, yet diverges with the emergence of a second branch above the initial solution branch. The resemblance ends there. Due to the heightened variation of the concavity of $\sin^3(t)$, additional solutions surface regularly, suggesting a perpetuation of this pattern, ultimately fostering a chaotic-like behavior in the quest for the Nash equilibrium.\\

Finally, due to the larger step size of $dt = 0.0002$ compared to $dt = 0.00005$, the precision of the solution endpoints is limited to approximately $a=56$. Despite the lack of accuracy beyond this threshold, $a$ values exceeding 56 still illustrate the evolution of the Nash equilibrium.\\

%% file: 9_FutureWork.tex
\section{Future Work}
We will try to understand the functions in the Nash equilibrium strategy more, and determine what is special about those particular functions. We will also examine optimal strategies and Nash equilibrium within other variations of the game, where we replace the $\sin(t)$ function with another odd function that isn't $\sin^3(t)$ (e.g., $\tan^{-1}(t)$), and dive deeper into why the functional and Nash equilibrium behave differently.\\

We also plan to try to optimize our code, whether that is from using a different programming language like C++, since Python is sluggish, or optimizing our code in other areas. Additionally, we want to employ superior approximation methods to reduce the amount of iterations necessary for reaching the Nash equilibrium and improve our accuracy.\\

Finally, let's delve into the diverse applications of this game theory problem. Economics stands out as a potential field of application due to the robust correlation between game theory and market dynamics. Another promising avenue is Artificial Intelligence (AI), given its focus on maximizing payoffs. A Nash equilibrium framework could offer insights into more subtle aspects. However, implementing this in AI applications requires a mechanism for translating various data types into a function and vice versa. Furthermore, while the potential for AI applications is possible, it's essential to acknowledge a significant limitation: this originated as a pure mathematical problem, and its broader applications still remain largely uncharted territory.

%% file: 10_Acknowledgements.tex
\section{Acknowledgements}
Special thanks to our mentors Dr. Hubert Bray (Duke University) and Dr. Dan Teague (NCSSM) for their guidance and support, as well as Dr. Michael Lavigne (NCSSM) for his assistance and advice. The authors would also like to thank NCSSM for giving them this research opportunity.